\documentclass[oneside, 12pt]{amsart}

\usepackage{amscd, amssymb, amsmath, mathrsfs}
\usepackage[english]{babel}
\usepackage[all,cmtip]{xy}
\usepackage[colorlinks, linktocpage, citecolor = blue, linkcolor = blue]{hyperref}
\usepackage{booktabs}
\usepackage{placeins}

\newcommand\mylabel[1]{\label{#1}}

\newtheorem{theorem}{Theorem}[section]
\newtheorem*{maintheorem}{Theorem}
\newtheorem{lemma}[theorem]{Lemma}
\newtheorem{proposition}[theorem]{Proposition}

\theoremstyle{definition}
\newtheorem{definition}[theorem]{Definition}

\newtheorem*{acknowledgement}{Acknowledgement}

\theoremstyle{remark}

\newcommand{\ZZ}	{\mathbb{Z}}
\newcommand{\QQ}	{\mathbb{Q}}
\newcommand{\RR}	{\mathbb{R}}

\newcommand{\PP}	{\mathbb{P}}

\newcommand{\ideala}    {\mathfrak{a}}

\newcommand  {\shN}     {\mathscr{N}}
\newcommand  {\shL}     {\mathscr{L}}

\newcommand{\Bl}	{\operatorname{Bl}}

\newcommand{\Cart}	{\operatorname{Cart}}

\newcommand{\Cl}	{\operatorname{Cl}}

\newcommand{\Exc}	{\operatorname{Exc}}

\newcommand{\lra}	{\longrightarrow}

\newcommand  {\Mat}     {\operatorname{Mat}}

\newcommand  {\maxid}   {\mathfrak{m}}

\renewcommand{\O}       {\mathscr{O}}

\newcommand  {\Pic}     {\operatorname{Pic}}

\newcommand  {\Proj}    {\operatorname{Proj}}

\newcommand  {\quadand} {\quad\text{and}\quad}

\newcommand  {\ra}      {\rightarrow}

\newcommand  {\Spec}    {\operatorname{Spec}}

\newcommand  {\Supp}    {\operatorname{Supp}}

\newcommand{\mydate}{
\number\day\space
\ifcase\month \or January\or February\or March\or April\or May\or June\or July\or August\or September\or October\or November\or December\fi 
\space\number\year
}

\begin{document}

\title[Mumford's rational pullback for Weil divisors]
      {A higher-dimensional generalization of Mumford's rational pullback for Weil divisors}

\author[Stefan Schr\"oer]{Stefan Schr\"oer}
\address{Mathematisches Institut, Heinrich-Heine-Universit\"at,
40204 D\"usseldorf, Germany}
\curraddr{}
\email{schroeer@math.uni-duesseldorf.de}


\dedicatory{17 July 2018}

\begin{abstract}
Mumford defined a rational pullback for Weil divisors
on normal surfaces, which  is linear,
respects effectivity, and satisfies the projection formula.
In higher dimensions, the existence of small resolutions of singularities precludes
such  general results.
We single out a higher-dimensional situation that resembles the
surface case and show for it that a rational pullback for Weil divisors   exists,
which is also linear, respects effectivity, and satisfies the projection formula.
\end{abstract}

\maketitle
\tableofcontents

%===========================================================
\section*{Introduction}
\mylabel{}

An important result for the theory of algebraic surfaces
is \emph{Mumford's rational pullback} \cite{Mumford 1961}: Let $f:X\ra S$ is the resolution of singularities
for a normal surface $S$. Then the   pullback for Cartier divisors extends to a unique map 
$$
f^*:Z^1(S)\lra Z^1(X)_\QQ=Z^1(X)\otimes_\ZZ\QQ
$$
for Weil divisors that retains the usual properties,
namely the map $f^*$ is linear, respects effectivity, and satisfies the projection formula.
The existence of this pullback relies on the fact that the intersection matrix $\Phi=(E_i\cdot E_j)$
for the exceptional divisors $E_i\subset X$ is negative-definite, 
with strictly negative   entries along the diagonal and positive off-diagonal entries.
Many definitions and results for smooth surfaces extend to normal surfaces, by using Mumford's rational pullback.
For example, the $\ZZ$-valued intersection form for Cartier divisors on proper surfaces
extends to a $\QQ$-valued intersection form for Weil divisors.
 
Despite the importance of a rational pullback, in particular for canonical divisors in the minimal model program,
there have been little attempts to generalize  Mumford's rational pullback to higher dimensions.
Indeed, the existence of \emph{small resolutions} of singularities in dimension $d\geq 3$ precludes
unconditional results. Nevertheless, de Fernex and Hacon \cite{de Fernex; Hacon 2009}
succeeded to construct a \emph{real-valued pullback} using valuation theory and asymptotic behavior in a surprising way. In this general set-up,
however, it is not so clear when linearity holds, effectivity is preserved
and the projection formula remains true.
The main goal of this paper is to single out a higher-dimensional situation
that  sufficiently resembles the surface case  and that yields a rational pullback with   these three   properties. 

We work in the following general set-up: Let $S$ be the spectrum of a   noetherian ring $R$ that
is local and  normal  of dimension $d\geq 2$, and $f:X\ra S$ be a proper birational morphism with $X$ integral and normal.
We do not require a ground field, but for the sake of exposition  we   assume that $R$ is excellent. Our main result is:

\begin{maintheorem}
{ (See Thm.\ \ref{criterion rational pullback})}
The morphism $f:X\ra S$ admits a rational pullback provided the following three conditions holds:
\begin{enumerate}
\item
All local rings $\O_{X,x}$ are $\QQ$-factorial.
\item
The exceptional locus $\Exc(X/R)$ and the closed fiber $f^{-1}(z)$ coincide  as closed sets,
and this is equidimensional of dimension $d-1$.
\item
Its irreducible components $E_1,\ldots, E_r$ have Picard number $\rho =1$.
\end{enumerate}
\end{maintheorem}

Note that the three assumptions hold in dimension $d=2$ for any resolution of singularities.
It is easy to produce example in higher dimension, by contracting suitable Cartier divisors.

The  key idea for the above result is to  work with the \emph{non-symmetric} square matrix $\Phi=(E_i\cdot C_j)$, where $C_j\subset E_j$ are chosen curves.
The crucial point is to establish that $A=-{}^t\Phi$ is an \emph{invertible M-matrix}, a very useful notion 
from linear algebra going back to Minkowski that generalizes positive-definiteness for symmetric matrices to
arbitrary square matrices. 
The theory of invertible M-matrices is widespread in applied mathematics, but perhaps not so well-known in pure mathematics.
Our approach relies on some recent contractibility results in \cite{Schroeer 2017}, which in turn are based on Cutkosky's 
study of \emph{graded linear system} and his generalization  of \emph{big invertible sheaves} to  non-integral schemes \cite{Cutkosky 2014}.

\medskip
The paper is organized as follows:
In the first section, we recall Mumford's rational pullback for surfaces, discuss
the problem of extending it to higher dimensions, and state our main result.
The second section contains the  proofs.

\begin{acknowledgement}
This research was conducted in the framework of the   research training group
\emph{GRK 2240: Algebro-geometric Methods in Algebra, Arithmetic and Topology}, which is funded
by the DFG. 
\end{acknowledgement}

%===========================================================
\section{Rational pullback}
\mylabel{Rational pullback}

Let $R$ be a local noetherian ring that is normal of dimension $d\geq 2$, with residue field $k=R/\maxid_R$,  spectrum $S$ and  closed
point $z\in S$. For the sake of exposition, I also assume that the ring $R$ is excellent.
Let $Z^1(S)$ be the group of  Weil Divisors, which is the free abelian group
generated by the prime divisors $D\subset S$.
This is an  \emph{ordered group} in the sense of \cite{A 4-7}, Chapter VI, where the  
\emph{positive elements} $D\geq 0$ are the effective divisors $D\subset X$.
The same applies to the group of $\QQ$-divisors $Z^1(S)_\QQ=Z^1(S)\otimes\QQ$ and the subgroup of Cartier divisors
$\Cart(S)$. Note that throughout we use the term ``positive'' in Bourbaki's sense $x\geq 0$,
and   ``strictly positive'' for $x>0$.

Let $f:X\ra S$ be a proper birational morphism, where $X$ is integral and normal.
For each point $s\in S$, the fiber $f^{-1}(s)$ is a proper scheme over the field $\kappa(s)$.
Its one-dimensional closed subschemes   are called \emph{vertical curves}.
For each invertible sheaf $\shL\in\Pic(X)$ and each vertical curve $C\subset f^{-1}(s)$, $s\in S$
we get an intersection number $(\shL\cdot C)=\chi(\shL_C)-\chi(\O_C)$,
where the Euler characteristics are computed over the residue field $\kappa(s)$.

For each effective Cartier divisor $D\subset S$, the subscheme  $f^{-1}(D)\subset X$ remains an effective Cartier divisor,
because $X$, $S$ are integral and $f:X\ra S$ is dominant.
This fact yields a pullback homomorphism  for Cartier divisors 
\begin{equation}
\label{pullback cartier}
f^*:\Cart(S)\lra \Cart(X)\subset Z^1(X),
\end{equation}
which is \emph{linear, increasing, and satisfies the  projection formula}. The latter means
$(\shL\cdot C) =0$ for each vertical curve $C\subset f^{-1}(s)$, $s\in S$.
We seek to extend \eqref{pullback cartier} to a \emph{rational pullback}
$f^*:Z^1(S)\ra Z^1(X)_\QQ$ that is also linear, increasing and satisfies the projection formula. 
Such an extension   exists a priori on the subgroup of $\QQ$-Cartier divisor.
The crux here is that we want to extend further,  without making any assumption on the \emph{class group} $\Cl(R)=Z^1(S)/\Cart(S)$.
However, to make sense of the intersection numbers in
the projection formula, we will usually assume that the local rings $\O_{X,x}$ are $\QQ$-factorial, that is, the
abelian groups $\Cl(\O_{X,x})$ are torsion groups.
 
Now suppose we are in dimension $d=2$. Then the \emph{exceptional locus} $\Exc(X/R)=\Supp(\Omega^1_{X/R})$ coincides with the closed fiber $f^{-1}(z)$,
and the underlying reduced closed subscheme $E\subset X$ is equidimensional, of dimension $\dim(E)=1$.
Decompose $E=E_1+\ldots+E_r$ into irreducible components. Under the assumption that all local rings $\O_{X,x}$ are $\QQ$-factorial,
we get a $\QQ$-valued intersection matrix $\Phi=(E_i\cdot E_j)$.
This   matrix is symmetric and negative-definite, an observation going back to 
Mumford \cite{Mumford 1961}, Artin \cite{Artin 1962} and Deligne (\cite{SGA 7b}, Expos\'e X, Corollary 1.9),
in various forms of generality.
Mumford used this to define the rational pullback $f^*:Z^1(S)\ra Z^1(X)_\QQ$ as follows:
For each prime divisor $D\subset S$, the strict transform $D'\subset X$
yields certain intersection numbers $(D'\cdot E_i)\geq 0$. Since $\Phi$ is invertible,
there is are unique rational numbers $m_1/n,\ldots,m_r/n$ with $(nD'\cdot E_j)= -(\sum m_i E_i\cdot E_j)$,
for each $1\leq j\leq r$. Mumford sets
$$
f^*(D) = D' + \frac{1}{n}\sum m_iE_i
$$
and extends by linearity (\cite{Mumford 1961}, Section II (b)). By construction, the projection formula $(f^*(D)\cdot E_j)=0$ holds.
A non-trivial fact from linear algebra ensures that all entries of $\Phi^{-1}$
are negative, hence $m_i/n\geq 0$, so the rational pullback   preserves effectivity.
 
In the general situation $d\geq 2$,  
and write $E_1,\ldots, E_r\subset \Exc(X/R)$ for the irreducible components of dimension $\dim(E_i)=d-1$.
This are precisely those prime divisors on $X$ whose images on $S$ cease to be a divisor.
The existence of a rational pullback can be seen as a problem in  \emph{linear programming}:

\begin{definition}
\mylabel{admits rational pullback}
Suppose that all local rings $\O_{X,x}$ are $\QQ$-factorial.
We say that  $f:X\ra S$ \emph{admits a rational pullback} if  
for each  $\shL\in\Pic(X)$ having a global section that does not vanish along $E_1,\ldots,E_r$,
there are unique rational number
$m_1/n,\ldots, m_r/n\in\QQ_{\geq 0}$ so that 
$
(\shL^{\otimes n}\cdot C) = -(\sum m_iE_i\cdot C)
$
for all    vertical curves $C\subset f^{-1}(s)$, $s\in S$.
\end{definition}

\noindent
Indeed, we then define the homomorphism $f^*:Z^1(S)\ra Z^1(X)_\QQ$ by the formula
\begin{equation}
\label{formula pullback}
f^*(D) = D' + \frac{1}{nn'}\sum m_iE_i,
\end{equation}
where $D$ is a prime divisor, $D'$ is its strict transform, $n'>0$ is an integer such that $n'D'$ is Cartier,
and the coefficients $m_i/n$ arise from
the invertible sheaf $\shL=\O_X(n'D')$. 
Clearly, this  does not depend on the choice of $n'$, and the map is linear, increasing, and satisfies the projection formula.
If $D\subset S$ is an effective Cartier divisor, then $f^{-1}(D)=D'+\sum m_iE_i$ is numerically trivial on all
vertical curves. Hence $f^{-1}(n'D)$ coincides with the rational pullback $f^*(n'D)$, whenver $n'D'\subset X$ is Cartier. 
Since $Z^1(X)_\QQ$ is torsion-free, we already have $f^{-1}(D)=f^*(D)$.
It follows that the rational pullback extends the usual pullback for Cartier divisors.
Conversely, if such a map $f^*:Z^1(S)\ra Z^1(X)_\QQ$ exists,    \eqref{formula pullback}  yields the desired
coefficients in Definition \ref{admits rational pullback}, by setting $\shL=\O_X(n'D')$.

It seems difficult  to verify directly that a morphism $f:X\ra S$ admits a rational pullback.
The main result of this paper is the following criterion, whose proof will occupy the
second section:

\begin{theorem}
\mylabel{criterion rational pullback}
The morphism $f:X\ra S$ admits a rational pullback provided the following three conditions holds:
\begin{enumerate}
\item
All local rings $\O_{X,x}$ are $\QQ$-factorial.
\item
The exceptional locus $\Exc(X/R)$ and the closed fiber $f^{-1}(z)$ coincide  as closed sets,
and this  is equidimensional of dimension $d-1$.
\item
Its irreducible components $E_1,\ldots, E_r$ have Picard number $\rho =1$.
\end{enumerate}
\end{theorem}

Each exceptional divisor $Y=E_i$ is a  proper $k$-scheme. Let $Z_1(Y)$ be the free abelian group generated by the integral
curves $C\subset Y$, and $\Pic(Y)\times Z_1(Y)\ra \ZZ$ be the ensuing intersection pairing.
The radical on the left is $\Pic^\tau(Y)$, the group of numerically trivial invertible sheaves.
By Finiteness of the Base, the residue class group $N^1(Y)=\Pic(Y)/\Pic^\tau(Y)$ is finitely generated.
Being torsion-free, it must be  free. Its rank $\rho\geq 0$ is called the \emph{Picard number}.
Let $ N^1(Y)\times N_1(Y)\ra\ZZ$ be the induced non-degenerate pairing.
Then also $N_1(Y)$ is finitely generated and free, of  rank $\rho\geq 0$.

In dimension $d=2$, we have $N_1(E_i)=\ZZ$, and the   conditions of the theorem are  automatically satisfied for any resolution of
singularities $f:X\ra S$. We thus recover  Mumford's   rational pullback.
Actually, it suffices to assume that the local rings $\O_{X,x}$ are $\QQ$-factorial.

It is not difficult to construct proper birational morphisms $f:X\ra S$ in arbitrary dimension $d\geq 2$
for which our result applies: Let $A$ be an excellent   discrete valuation ring, with  residue field
$k=A/\maxid_A$, and consider any projective flat $A$-scheme $Y$ whose closed fiber 
$Y\otimes_Ak$ is smooth, with Picard number $\rho=1$. 
According to \cite{Schroeer 2017}, Proposition 1.6, there is an effective Cartier divisor $Z\subset Y\otimes_Ak$
so that on the blowing-up $\varphi:\Bl_Z(Y)\ra Y$, the strict transform $E$ of the closed fiber $Y\otimes_Ak$
admits a contraction $\varphi':\Bl_Z(Y)\ra Y'$ to some projective $A$-scheme $Y'$.
Such a construction resembles the elementary transformations for projective bundles,
and was already used in \cite{Schroeer 1999} for surfaces.
Using Bertini, one can arrange things that $Z$ is smooth, hence the total space $\Bl_Z(Y)$ is regular.
The image $z=\varphi'(E)$ is a closed point. Let $R=\O_{Y',z}$ be the resulting local ring. With
$S=\Spec(R)$, the resulting base-change $X=\Bl_{Z}(Y)\times_{Y'}S$ yields a proper birational
morphism $f:X\ra S$ for which Theorem \ref{criterion rational pullback} applies, and thus admits a rational pullback
$f^*:Z^1(S)\ra Z^1(X)_\QQ$.

Let me close this section with a standard example that shows that in dimension $d\geq 3$, there are  many important  $f:X\ra S$ 
that  do not admit a rational pullback: Fix a ground field $k$, and let $R=k[[x,y,u,v]]/(xy-uv)$.
By taking partial derivatives, one sees that $S=\Spec(R)$ has an isolated singularity.
The ideal $\ideala=(x,u)$ defines a prime divisor $D\subset S$.
On the  blowing-up $X=\Proj(R[\ideala T])$ the Cartier divisor defined by the canonical inclusion $\O_X(1)\subset\O_X$
is the blowing-up of the spectrum of $R/\ideala=k[[y,v]]$ at the origin, as one sees by computing $D_+(xT)$ and $D_+(uT)$.
In turn, the exceptional locus $ \Exc(X/R)$ coincides with the closed fiber $f^{-1}(z)$, and is  is a copy $C=\PP^1$
of the projective line. In particular, it contains no Cartier divisor. Hence   $f:X\ra S$ is
a \emph{small resolution of singularities}. The invertible sheaf $\shL=\O_X(1)$ for the blowing-up is relatively ample, so $(\shL\cdot C)>0$.
Since $r=0$, the condition in Definition \ref{admits rational pullback} does not hold.

%===========================================================
\section{Invertible M-matrices}
\mylabel{Mumford's pullback}

After some preparation, we  now give a proof for Theorem \ref{criterion rational pullback}.
Notation is as in the previous section, in particular $f:X\ra S=\Spec(R)$ is a proper birational morphism,
where $R$ is a local noetherian ring that is normal and of dimension $d\geq 2$, and the scheme $X$ is integral and normal.
Let us start with a useful general observation:

\begin{lemma}
\mylabel{effective cartier}
There is an  effective Cartier divisors $D\subset S$ whose strict transform $D'\subset X$
intersects each  of the exceptional divisors $E_1,\ldots,E_r$.
\end{lemma}

\proof
Let $\eta_i\in E_i$ be the generic points. Then the local rings $\O_{X,\eta_i}$ are one-dimensional.
According to \cite{Gross 2012}, Corollary 1.6, there is a common affine open neighborhood $U\subset X$ for the points
$\eta_1,\ldots,\eta_r\in X$. The complement $A=X\smallsetminus U$, endowed with the reduced scheme structure,
is an effective Weil divisor that intersects all $E_i$ (confer \cite{Hartshorne 1966}, Chapter II, Proposition 3.1).
Its image $f(A)\subset S$ is an effective Weil divisor,
with strict transform $A$. Let $\ideala\subset R$ be the resulting ideal, and choose some non-zero $\zeta\in\ideala$.
This defines an effective Cartier divisor $D\subset S$ containing $f(A)$. In turn, the
strict transform $D'\subset X$ contains $A$, thus intersects each exceptional divisor $E_i$.
\qed 

\medskip
From now on, we suppose that conditions (i)--(iii) from Theorem \ref{criterion rational pullback} hold.
In particular, each Weil divisor on $X$ is $\QQ$-Cartier.
Moreover,   $N_1(E_i)\simeq \ZZ$. Actually, there is a canonical identification:

\begin{proposition}
\mylabel{ratio positive}
For each curve $C\subset E_i$,
the class $[C]\in N_1(E_i)$ is non-zero, and for each further curve  $C'\subset E_i$,
the     equation
$[C]=\mu[C']$ defines a ratio $\mu\in\QQ_{>0}$.
\end{proposition}

\proof
For this, it suffices to treat the case that both curves $C,C'$ are irreducible.
Since the proper $k$-scheme $E_i$ is connected, there is a sequence of irreducible curves 
$$
C=C_0, C_1,\ldots, C_n=C'
$$
with $C_{i}\cap C_{i+1}$ non-empty.
By induction on $n\geq 0$,  it suffices to treat the case that $C\cap C'$ is non-empty.
Choose an affine open neighborhood $U\subset X$ of some intersection point $a\in C\cap C'$.
The complement $D=X\smallsetminus U$ is an effective Weil divisor, and the intersections $D\cap C$ and $D\cap C'$
are both zero-dimensional, hence the intersection numbers $(D\cdot C)$ and $(D\cdot C')$ are strictly positive. 
It follows that both classes $[C], [C']$ are non-zero, and that the  ratio $\mu\in\QQ$ is strictly positive.
\qed

\medskip
In turn, the one-dimensional vector spaces $N_1(E_i)_\QQ$ are ordered groups  whose positive
elements are the   $\mu[C]$ with $\mu\geq 0$, where $C\subset E_i$ is any curve.
For each $1\leq i\leq r$, we now choose some curve $C_i\subset E_i$, 
and consider the resulting \emph{intersection matrix} 
$$
\Phi=(E_i\cdot C_j)\in\Mat_{r\times r}(\QQ).
$$
Note that this matrix is usually \emph{not symmetric}, in contrast to the situation in dimension $d=2$. Furthermore, it \emph{depends on the choices}
of curves.
It is easy to determine the signs in the intersection matrix, which are actually independent from the chosen curves:

\begin{lemma}
\mylabel{intersection matrix}
We have $(E_j\cdot C_j)<0$ for all $1\leq j\leq r$, and $(E_i\cdot C_j)\geq 0$ for $i\neq j$,
with equality if and and only if $E_i\cap E_j=\varnothing$.
\end{lemma}

\proof
If $E_i\cap E_j=\varnothing$, the intersection $E_i\cap C_j$ is also empty, hence $(E_i\cdot C_j)=0$.
In the   case  where $i\neq j$ and $E_i\cap E_j\neq\varnothing$,
we may choose the curve $C_j\subset X_j$ so that it is not contained in $E_i$ but intersects $E_i$.
In turn, we have $\dim(E_i\cap C_j)=0$, and thus $(E_i\cdot C_j)>0$.
It remains to show $(E_j\cdot C_j)<0$. 
Choose some non-zero non-unit $\zeta\in R$, with Cartier divisor $D\subset S$,
and decompose $f^*(D)= D' +\sum m_iE_i$, where $D'$ is the strict transform
and the coefficients are $m_i>0$. The inclusion $D'\cap E_j\subset E_j$  is strict, so we may choose
the curve $C_j\subset E_j$ not contained in $D'$.
Since $f^*(D)$ becomes numerically trivial on $C_j$, we have 
$$
-m_j(E_j\cdot C_j) = (D'\cdot C_j) + \sum_{i\neq j} m_i(E_i\cdot C_j).
$$
The intersection numbers on the right are positive, and we are done if at least one is strictly positive.
If $r>1$, we find some $i\neq j$ with $E_i\cap E_j\neq\varnothing$, because $f:X\ra S$ has connected fibers.
We may choose $C_j\subset E_j$ so that it   intersects $E_i$ but is not contained in $E_i$. Thus   $(E_i\cdot C_j)$ are strictly positive.
If $r=1$ we have $j=1$, and the intersection $D'\cap E_1$ is non-empty.
Now we choose  $C_1\subset E_1$ so that it intersects $D'$ but is not contained in $D'$. 
In both cases, one intersection number on the right is strictly positive.
\qed

\medskip
Given any real $r\times r$-matrix $A=(\alpha_{ij})$ whose off-diagonal entries are $\alpha_{ij}\leq 0$, we may write it  in the form $A=sE - B$,  
for some scalar $s>0$ and some matrix $B=(\beta_{ij})$ all whose entries are $\beta_{ij}\geq 0$. Here $E$ denotes the unit matrix.
Recall that the \emph{spectral radius} $\rho(B)\geq 0$ is the maximal   length occurring for the complex eigenvalues of $B$.
If     $\rho(B)<s$ for some scalar  $s>0$ ,
the matrix $A$   called an \emph{invertible M-matrix}.  Note that if this holds for some $s>0$, it also holds for
all $s'\geq s$.

The  terminology seems to refer to Minkowski, and such matrices have amazing properties.
Berman and Plemmons give fifty   characterizations of invertible M-matrices (\cite{Berman; Plemmons 1994}, Chapter 6, Theorem 2.3).
One of them is condition $(I_{27})$: There is a column vector ${}^t(x_1,\ldots,x_r)\in\RR^r_{>0}$ with $Ax\in\RR_{>0}^r$.
Another one is  $(N_{38})$: The inverse $A^{-1}=(\lambda_{ij})$ has all entries $\lambda_{ij}\geq 0$.
Note that for symmetric matrices, the notion boils down to \emph{positive-definiteness}.
We now apply this  theory  to the negative transpose of  our intersection matrices:

\begin{proposition}
\mylabel{inverse matrix}
The matrix $A=-{}^t\Phi$ is an invertible M-matrix. In particular, we have 
$\det(\Phi)\neq 0$, and all entries of the inverse matrix $\Phi^{-1}$
are negative.
\end{proposition}

\proof
Suppose $\det(\Phi)=0$. Then there is some non-zero Cartier divisor $\sum m_iE_i$ such
that the resulting invertible sheaf $\shN=\O_X(\sum m_iE_i)$ is numerically trivial on the proper $k$-scheme $E=E_1\cup\ldots\cup E_r$.
After passing to some multiple and renumeration, we may assume that the summands $m_iE_i$ are Cartier,
and that the non-zero coefficients are $m_1,\ldots, m_a>0$ and $m_b,\ldots, m_r<0$,
for some $1\leq a< b\leq r+1$. This gives effective Cartier divisors
$$
E'=\sum_{i=1}^am_iE_i \quadand E''=\sum_{i=b}^r(-m_i)E_i 
$$
whose invertible sheaves 
become numerically equivalent on   $E$. We have  $(E'\cdot C_i)\geq 0$ for all $i>a$, and  $(E''\cdot C_i)\geq 0$ for all $i<b$.  
Setting $Y=E'$, we see that the invertible sheaf $\shL=\O_X(Y)$ is nef on each $E_i$.
On the other hand, the restriction $\shL^\vee|Y=\O_Y(-Y)$ is a \emph{big invertible sheaf}, according to \cite{Schroeer 2017}, Theorem 1.5.
This means that the   homogeneous  spectrum of the graded ring 
$$
R(Y,\shL^\vee|Y) = \bigoplus_{n\geq 0} H^0(Y,\shL^{\otimes -n}_Y)
$$
attains the maximal possible dimension $\dim(Y)=d-1$. The notion of big invertible sheaves on integral schemes
is common. However, here it is crucial to work with Cutkosky's generalization \cite{Cutkosky 2014} to arbitrary proper schemes,
because our scheme $Y$ usually is reducible and non-reduced.
Also note that in loc.\ cit.\ we worked with  schemes that are proper over an excellent discrete valuation ring,
but the argument literally hold true over our excellent local ring $R$.

By \cite{Cutkosky 2014}, Lemma 10.1 combined with Lemma 9.1, there is some irreducible component $E_j\subset Y$ such that  $\shL^\vee|{E_j}$ is 
big. In particular, there is some integer $n>0$ and some non-zero global section $\sigma\in H^0(E_j,\shL^{\otimes -n}|E_j)$.
Write $Z\subset E_j$ for the resulting zero-locus, and choose an irreducible curve $C_j'\subset E_j$
not contained in $Z$ but intersecting $Z$. It follows that $(\shL^{\otimes -n}\cdot C_j')>0$,
contradicting that $\shL|E_j$ is nef.
Thus $\det(\Phi)\neq 0$.

To understand  the inverse matrix $A^{-1}= -{}^t\Phi^{-1}$, choose an effective Cartier  divisor $D\subset S$
as in Lemma \ref{effective cartier}, and write $f^*(D)=D'+\sum m_iE_i$. Here $D'$ is the strict transform,
and all coefficients $m_i$ and intersection numbers $\lambda_j=(D'\cdot C_j)$ 
are strictly positive. Moreover,  $- (\sum m_iE_i\cdot C_j) = (D'\cdot C_j)$ for each $1\leq j\leq r$.
In terms of matrix multiplication, this means
\begin{equation}
\label{semipositive}
(m_1,\ldots, m_r) \cdot (-\Phi) = (\lambda_1,\ldots,\lambda_r).
\end{equation}
In turn, $A=-{}^t\Phi$ sends
the transpose of  $(m_1,\ldots,m_r)$ to the transpose $(\lambda_1,\ldots,\lambda_r)$, and 
all entries of these vectors are strictly positive.
According to  $(I_{27})$
in \cite{Berman; Plemmons 1994}, Chapter 6, Theorem 2.3 our $A$ is an invertible  M-matrix, 
and this ensures by $(N_{38})$ that the entries in $A^{-1}=-{}^t\Phi^{-1}$ are positive.
\qed

\medskip
\emph{Proof of Theorem \ref{criterion rational pullback}.}
Let $\shL$ be an invertible sheaf on $X$ having a global section that does not
vanish on any exceptional divisor $E_1,\ldots, E_r$. The corresponding Cartier divisor $D'\subset X$
is the strict transform of the Weil divisor $D=f(D')$.
The inclusions $D'\cap E_i\subset E_i$ are strict, so we may choose the curves $C_i\subset E_i$
so that they are not contained in $D'$. In turn, we have $\lambda_j=(D'\cdot C_j)\geq 0$.
Condition (i) ensures that each Weil divisor $E_i\subset X$ is $\QQ$-Cartier, so we  may 
form the intersection matrix $\Phi=(E_i\cdot C_j)$. By Proposition \ref{inverse matrix},
this matrix is invertible, and the entries of its inverse are negative.
The equation
$$
(m_1/n,\ldots, m_r/n) = (\lambda_1,\ldots,\lambda_r)\cdot (-\Phi^{-1})
$$
defines rational numbers $m_i/n\in\QQ_{\geq 0}$.
In turn, we have
\begin{equation}
\label{intersection correct}
(\shL^{\otimes n}\cdot C_j) = (nD'\cdot C_j) =n\lambda_j=   -\sum m_i(E_i\cdot C_j) = -(\sum m_iE_i\cdot C_j),
\end{equation}
and the $m_i/n\in\QQ_{\geq 0}$ are the only fractions having this property.
By condition (ii), any vertical curve  $C\subset f^{-1}(s)$ lies over the closed point $s=z$.
If $C$ is contained in $E_j$, we have $[C]=\mu_j[C_j]$ inside $N_1(E_j)$ for some ratio $\mu_j>0$, according Proposition \ref{ratio positive}.
Note that the latter  relies on condition (iii).
From equation \eqref{intersection correct}, we infer $(\shL^{\otimes n}\cdot C)=-(\sum m_iE_i\cdot C)$.
Thus the condition  in Definition \ref{admits rational pullback}  is fulfilled, in other words,  our morphism $f:X\ra S$ admits a rational pullback.
\qed

%===========================================================

\end{document}